\documentclass{amsart}
\usepackage{amssymb,euscript,amsmath, mathrsfs, amscd, mathabx}
\usepackage[dvips]{graphicx}
\usepackage[dvips]{color}
\usepackage{bm}

\newcounter{ENUM}
\newcommand{\itm}{\item}

\newenvironment{Ilist}{\renewcommand{\theENUM}{\Roman{ENUM}}\renewcommand{\itm}{\addtocounter{ENUM}{1}\item[(\theENUM)]}\begin{itemize}\setcounter{ENUM}{0}}{\end{itemize}}

\def\risom{\overset{\sim}{\rightarrow}}

\input xy
\xyoption{all}
\CompileMatrices

\def\ZZ{{\mathbb Z}}

\def\PP{{\mathbb P}}

\def\bn{{\bm{n}}}

\def\sL{{\mathscr L}}

\def\sO{{\mathscr O}}

\def\fg{{\mathfrak g}}

\def\e{\varepsilon}
\def\vp{\varphi}

\def\ch{\operatorname{char}}

\def\id{\operatorname{id}}

\def\ord{\operatorname{ord}}

\def\lcm{\operatorname{lcm}}

\newtheorem{thm}{Theorem}[section]
\newtheorem{prop}[thm]{Proposition}
\newtheorem{lem}[thm]{Lemma}
\newtheorem{cor}[thm]{Corollary}

\theoremstyle{definition}
\newtheorem{defn}[thm]{Definition}

\newtheorem{ex}[thm]{Example}
\newtheorem{sit}[thm]{Situation}
\theoremstyle{remark}
\newtheorem{notn}[thm]{Notation}
\newtheorem{rem}[thm]{Remark}

\numberwithin{equation}{section}
\numberwithin{figure}{section}

\begin{document}
\title[Dimension counts for limit linear series]{Dimension counts for limit 
linear series on curves not of compact type}
\author{Brian Osserman}
\begin{abstract} 
We first prove a generalized Brill-Noether theorem for linear series with 
prescribed multivanishing sequences on smooth curves. We then apply this 
theorem to prove that spaces of limit linear series have the expected 
dimension for a certain class of curves not of compact type, whenever the 
gluing conditions in the definition of limit linear series impose the 
maximal codimension. Finally, we investigate these gluing conditions in
specific families of curves, showing expected dimension in several 
cases, each with different behavior. One of these families sheds new
light on the work of Cools, Draisma, Payne and Robeva in tropical
Brill-Noether theory.
\end{abstract}

\thanks{The author was partially supported by NSA grant H98230-11-1-0159
and Simons Foundation grant \#245939 during the preparation of this work.}
\maketitle

\section{Introduction}

In \cite{os25}, the author introduced a theory of limit linear series
for nodal curves not of compact type, along with an equivalent 
definition, generalizing the Eisenbud-Harris definition, for curves
of `pseudocompact type', meaning that their dual graphs yield trees
after collapsing all multiple edges. It was shown that linear series
always specialize to limit linear series, and, as in the work of Eisenbud and
Harris, for curves of pseudocompact type it was shown that limit linear 
series occurring in families of the expected dimension can be smoothed
to linear series on smooth curves. In order to show that the theory is
useful, it thus remains to show first that it is tractable, and 
second that -- at least in some interesting cases -- one does in fact
obtain families of limit linear series having the expected dimension.
This is precisely what is accomplished in the present paper.

We work in the context of curves of pseudocompact type, where a basic
ingredient in the definition of limit linear series is the notion of
`multivanishing sequence', defined as follows:

\begin{defn}\label{def:multivanishing} Let $X$ be a smooth projective
curve, $r,d \geq 0$, and $D_0 \leq D_1 \leq \dots \leq D_{b+1}$ a 
sequence of effective divisors on $X$, with $D_0=0$ and $\deg D_{b+1}>d$. 
Given $(\sL,V)$ a $\fg^r_d$ on $X$, define the
\textbf{multivanishing sequence} of $(\sL,V)$ along $D_{\bullet}$ to
be the sequence
$$a_0 \leq \dots \leq a_r$$
where a value $a$ appears in the sequence $m$ times if for some $i$ we
have $\deg D_i=a$, $\deg D_{i+1}>a$, and 
$\dim \left(V(-D_i)/V(-D_{i+1})\right)=m$.
\end{defn}

In the above, if $D$ is an effective divisor and $(\sL,V)$ a linear series,
we write $V(-D)$ for $V \cap \Gamma(X,\sL(-D))$.

Thus, this generalizes usual vanishing sequences, but also addresses a
wide variety of geometric conditions, from secancy (in the sense of
requiring two or more points to map to the same point) to bitangency,
and so forth. The definition of limit linear series for curves of
pseudocompact type (recalled in Definition \ref{defn:lls} below) involves 
a condition on multivanishing sequences, as well as a gluing condition.

Our first result, given precisely in Theorem \ref{thm:bn-genl}
below, is a generalized Brill-Noether theorem for linear series with
imposed multivanishing sequences. We then apply this theorem to 
investigate which reducible nodal curves are \textbf{Brill-Noether
general} -- that is, have limit linear series spaces of the expected
dimension. We prove in Theorem \ref{thm:exp-dim} that any curve of 
pseudocompact type whose components are Brill-Noether general with respect 
to imposed multivanishing is itself Brill-Noether general 
as long as the gluing conditions impose maximal
codimension. These constitute the general results of the paper.

The remainder of the paper consists of studying two extremes where gluing
becomes relatively tractable: when there are few nodes (at most two or 
three) between any given pair of components, or at the opposite extreme,
`binary curves' consisting of two rational components glued to one another
at $g+1$ nodes. In each situation, we are able to prove the necessary
independence of gluing conditions under suitable hypotheses. Specifically,
in order to define limit linear series on a given curve, we first choose
two additional structures: a numerical ``chain structure'' which can be
thought of as describing the singularities of a one-parameter smoothing
of the curve, and an ``enriched structure'' consisting of a collection of
twisting bundles which may likewise be obtained from a (regular) smoothing.
Generality of enriched structure corresponds to a more typical 
algebrogeometric notion of generality. Imposing conditions on chain 
structures in some sense restricts to more special families of curves, but 
insofar as it can be used to restrict directions of approach to a given 
nodal curve, it can still be used to ensure that a space of limit linear
series has the expected dimension. Indeed, the notion of generality in
tropical Brill-Noether theory (as in Cools-Draisma-Payne-Robeva 
\cite{c-d-p-r}) is essentially the same as restricting chain structures
-- see the end of \S \ref{sec:review}.

In the case of few nodes, in Corollary \ref{cor:few-nodes} we produce 
families of curves which are Brill-Noether general
when restrictions are placed on the chain structures, 
irrespective of enriched structure.
This may be viewed as a generalization of \cite{c-d-p-r}, and in fact
helps explain why their genericity condition should lead to good behavior.
Then, in Corollary
\ref{cor:two-node-two} we show that a certain narrower family of curves
has the complementary behavior, with Brill-Noether generality
occurring for general enriched structures,
irrespective of the chain structure. Here, the family in question contains
in particular the curves considered by Jensen and Payne in their work 
\cite{j-p1} on a tropical approach to the Gieseker-Petri theorem. 
Finally, in Corollary \ref{cor:binary} we show that for binary curves,
the behavior is very different: whether or not the curve is Brill-Noether
general is independent of both enriched
structures and chain structures, and may be expressed in terms of usual
linear series (in a certain restricted range of multidegrees) on the
underlying curve. Using work of Caporaso \cite{ca3} we conclude that
limit linear series spaces for $r \leq 2$ have the expected dimension 
on general binary curves. 

This last result draws a distinction between our limit linear series
and the notion introduced by Amini and Baker in \cite{a-b1}: a general
binary curve of genus at least $3$ is not `hyperelliptic' with respect
to our definition of limit linear series, but it is with respect to the
Amini-Baker definition. In fact, for us the main purpose of the final
sections of the paper is to validate the definition of limit linear series 
given in \cite{os25} by producing quite distinct infinite families
of curves where the limit linear series spaces have the correct dimension.
At the same time, as mentioned above, the particular families studied have 
also arisen (at least in special cases) elsewhere in the literature, leading 
to natural points of contact. Indeed, the relationship between Corollary 
\ref{cor:few-nodes} and the work in \cite{c-d-p-r} is highly suggestive
on the tropical side, pointing to families of graphs which are and are not 
likely to be Brill-Noether general in the tropical sense. In addition, 
we hope that our work will lead to a
new proof that the tropical linear series studied in \cite{c-d-p-r} can
all be smoothed, as was recently proved by Cartwright, Jensen and Payne
in \cite{c-j-p1}. We will investigate this and other aspects of the
relationship to the tropical theory in \cite{os24}. Finally, our work leads 
to new explicit criteria for Brill-Noether generality in terms of special
fibers of degenerations, stated in Corollary \ref{cor:few-nodes-genl}.

\subsection*{Acknowledgements} I would like to thank Eduardo Esteves, Melody
Chan, Matt Baker and Sam Payne for helpful conversations.

\subsection*{Conventions}
We do not assume base fields to be algebraically closed, but we do 
assume that our nodal curves are split over their base field, meaning
that all components and nodes are defined over the base field.

If $v$ is the vertex of the dual graph of a nodal curve, we let $Z_v$
denote the corresponding component of the curve. If $e$ is an edge,
we let $P_e$ denote the corresponding node.

If $e$ is an edge of a directed graph, we denote by $h(e)$ and $t(e)$
the head and tail of $e$ respectively.

\section{Background on limit linear series}\label{sec:review}

In this section we recall the definition of limit linear series 
introduced in \cite{os25}, simplifying somewhat as a consequence of
restricting to the context of curves of pseudocompact type. Although
our presentation is mathematically self contained, we refer the reader
to \S 2 of \cite{os25} for additional remarks and examples. However, for 
the benefit of readers coming from
tropical geometry we do include a rough dictionary between our definitions 
and those arising in tropical Brill-Noether theory at the end of this 
section.

We begin with some definitions of a combinatorial nature. In the below,
$\Gamma$ will be obtained by choosing a directed structure on the dual
graph of a curve of pseudocompact type. 
The following definition forms the basis for our approach to keeping 
track of chains of rational curves inserted at the nodes of the original
curve.

\begin{defn} A \textbf{chain structure} on $\Gamma$ is a function 
$\bn:E(\Gamma) \to \ZZ_{> 0}$. 
\end{defn}

The chain structure will determine the length of the chain of rational
curves inserted at a given node; for notational convenience, the
trivial case (in which no rational curves are inserted) corresponds to
$\bn(e)=1$.

We work in the following situation throughout.

\begin{sit}\label{sit:gamma} Let $\Gamma$ be a directed graph without
loops, and $\bn$ a chain structure on $\Gamma$. For each pair of an edge 
$e$ and adjacent vertex $v$ of $\Gamma$, 
let $\sigma(e,v)=1$ if $e$ has tail $v$, and $-1$ if $e$ has head $v$.

Let $\bar{\Gamma}$ be the graph obtained from $\Gamma$ by collapsing
all multiple edges, and assume that $\bar{\Gamma}$ is a tree. 
\end{sit}

\begin{defn} An \textbf{admissible multidegree} $w$ of total 
degree $d$ on $(\Gamma,\bn)$ consists of a function 
$w_{\Gamma}: V(\Gamma) \to \ZZ$ together with a tuple
$(\mu(e))_{e \in E(\Gamma)}$, where each $\mu(e) \in \ZZ/\bn(e)\ZZ$,
such that
$$d = \#\{e \in E(\Gamma): \mu(e) \neq 0\}
+ \sum_{v \in V(\Gamma)} w_{\Gamma}(v).$$
\end{defn}

The idea behind admissible multidegrees is that in order to extend line
bundles, we need only consider multidegrees which have degree $0$ or $1$
on each rational curve inserted at the node, with degree $1$ occurring at
most once in each chain. Thus, $\mu(e)$ determines where on the chain
(if anywhere) positive degree occurs. 
See Definition \ref{def:multideg-lb} below for details.

We now define twists of multidegrees.

\begin{defn}\label{def:twist-node} If $(e,v)$ is a pair of an edge $e$ and 
an adjacent vertex $v$ of $\bar{\Gamma}$, 
given an admissible multidegree $w$, we define the \textbf{twist} of
$w$ at $(e,v)$ to be obtained from $w$ as follows:
for each $\tilde{e}$ of $\Gamma$ over $e$, increase $\mu(\tilde{e})$ by 
$\sigma(\tilde{e},v)$. Now, decrease $w_{\Gamma}(v)$ by the number of 
$\tilde{e}$ for which $\mu(\tilde{e})$ had been equal to $0$, and for each 
$\tilde{e}$, if the new $\mu(\tilde{e})$ is zero, increase 
$w_{\Gamma}(v')$ by $1$, where $v'$ is the other vertex adjacent to $v$.
\end{defn}

Twists will be the change in multidegrees accomplished by twisting by
certain natural line bundles; see Notation \ref{not:more-twist} below.

\begin{defn} An admissible multidegree $w$ is \textbf{concentrated} at
a vertex $v \in V(\Gamma)$ if for each $v' \neq v$, we have that $w$ 
is negative in index $v'$ after twisting at $(e,v')$, where $e$ is the
edge of $\bar{\Gamma}$ from $v'$ in the direction of $v$.
\end{defn}

For the sake of simplicity, the above definition is slightly more
restrictive than that of \cite{os25}; see Remark \ref{rem:concen-diff}.
It generalizes the situation for Eisenbud and Harris of considering
line bundles with degree $d$ on one component and degree $0$ on the others.
If an admissible multidegree $w_v$ is concentrated at $v$ and also has
negative degree in index $v$, then it will lead to a vacuous theory of
limit linear series, because of Proposition 3.3 of \cite{os25}. Accordingly,
we will always assume implicitly that any such $w_v$ is nonnegative in
index $v$.

We will work throughout in the following situation.

\begin{sit}\label{sit:concen} Suppose we are given 
an admissible multidegree $w_0$, and let 
$(w_v)_{v \in V(\Gamma)}=(w_v,\mu_v(\bullet))_v$ be a
collection of admissible multidegrees, each obtained by $w_0$ by successive
twists, and such that:
\begin{Ilist}
\itm each $w_v$ is concentrated at $v$;
\itm for each $v,v' \in V(\bar{\Gamma})$ connected by an edge $e$, the
multidegree $w_{v'}$ is obtained from $w_v$ by twisting $b_{v,v'}$ times at 
$(e,v)$, for some $b_{v,v'}\in \ZZ_{\geq 0}$.
\end{Ilist}
\end{sit}

The below graph will help us keep track of the relevant multidegrees.

\begin{defn}\label{def:tree-graph}
In Situation \ref{sit:concen}, let 
$$V(\bar{G}(w_0)) \subseteq \ZZ^{V(\Gamma)} \times \prod_{e \in E(\Gamma)} 
\ZZ/\bn(e)\ZZ$$
consist of admissible multidegrees $w$ such that there exist 
$v,v' \in V(\bar{\Gamma})$ connected by some edge $e$, with $w$ obtainable 
from $w_v$ by twisting $b$ times at $(e,v)$, for some $b$ with
$0 \leq b \leq b_{v,v'}$.

There is an edge $\epsilon$ from from $w$ to $w'$ in $\bar{G}(w_0)$ if there
exists $(e,v)$ in $\bar{\Gamma}$ such that $w'$ is obtained from $w$ by
twisting at $(e,v)$.
\end{defn}

Thus, $\bar{G}(w_0)$ is essentially a tree, obtained from $\bar{\Gamma}$ by 
subdividing each edge into $b_{v,v'}$ edges, and replacing each resulting 
edge with a pair of edges going in opposite directions.

We now move on to the geometric constructions which underlie our definition
of limit linear series. First, in the non-compact-type case, additional
structure beyond the underlying nodal curve is necessary in order to make
a useful definition.

\begin{defn}\label{def:enriched} If $X'$ is a projective nodal curve with 
dual graph $\Gamma'$, an
\textbf{enriched structure} on $X'$ consists of the data, for each 
$v \in V(\Gamma')$ of a line bundle $\sO_v$ on $X'$,
satisfying the following conditions:
\begin{Ilist}
\itm for any $v \in V(\Gamma')$, we have
$$\sO_v|_{Z_v} \cong \sO_{Z_v}(-(Z_v^c \cap Z_v)),\text{ and }
\sO_v|_{Z_v^c} \cong \sO_{Z_v^c}(Z_v^c \cap Z_v);$$
\itm we have
$$\bigotimes_{v \in V(\Gamma')} \sO_v \cong \sO_{X}.$$
\end{Ilist}
In the above, $Z_v^c$ is the closure of the complement of $Z_v$.
\end{defn}

The curve on which we place an enriched structure will not be the
original nodal curve, but the following curve obtained by also taking
the chain structure into account.

\begin{defn}
Given a projective nodal curve $X$ with dual graph $\Gamma$ and a chain 
structure $\bn$, let $\widetilde{X}$ denote 
the nodal curve obtained from $X$ by, for each $e \in E(\Gamma)$,
inserting a chain of $\bn(e)-1$ projective lines at the corresponding
node. Let $\widetilde{\Gamma}$ be the dual graph of $\widetilde{X}$,
with a natural inclusion $V(\Gamma) \subseteq V(\widetilde{\Gamma})$.
\end{defn}

Finally, our definition of limit linear series will occur in the following
context.

\begin{sit}\label{sit:geom} In Situation \ref{sit:concen}, suppose further
that we have a projective nodal curve $X$ over a field $k$ with dual graph 
$\Gamma$, and an enriched structure $(\sO_v)_v$ on the corresponding 
$\widetilde{X}$. Fix also $r>0$, and let $d$ be the total multidegree of 
$w_0$.
\end{sit}

We now describe how our combinatorial notions of multidegrees and twists
arise in the geometric setting.

\begin{defn}\label{def:multideg-lb} Using our orientation of $E(\Gamma)$,
an admissible multidegree $w$ of total degree $d$ on $(X,\bn)$ gives a 
multidegree of total degree $d$ on $\widetilde{X}$
by assigning, for each $e \in E(\Gamma)$, 
degree $0$ on each component of the corresponding chain of projective curves,
except for degree $1$ on the $\mu(e)$th component when $\mu(e) \neq 0$. 
\end{defn}

The following notation will not be used later, and is necessary only
to set up Notation \ref{not:more-twist}.
In Situation \ref{sit:geom}, for any edge $\e \in E(\bar{G}(w_0))$, starting 
at $w=(w_{\Gamma},(\mu(e))_{e \in E(\Gamma)})$ and determined by twisting at
$(e,v)$, let $\widetilde{\Gamma}'$ be the graph
obtained from $\widetilde{\Gamma}$ by removing, for each 
edge $\widetilde{e}$ of $\Gamma$ lying over $e$, the 
$(\sigma(\widetilde{e},v)\mu(\widetilde{e})+1)$st edge of 
$\widetilde{\Gamma}$ lying over $\widetilde{e}$, starting from $v$.
Then let $S\subseteq V(\widetilde{\Gamma})$ consist of the vertices in the 
connected component of $\widetilde{\Gamma}'$ containing $v$.
Next, let $\sO_{\e}$ be the twisting line bundle on $\widetilde{X}$ defined by
$$\sO_{\e}:=\bigotimes_{v' \in S} \sO_{v'}.$$
Similarly, given $w,w' \in V(\bar{G}(w_0))$, let $P=(\e_1,\dots,\e_m)$ be
a minimal path from $w$ to $w'$ in $\bar{G}(w_0)$, and set
$$\sO_{w,w'}=\bigotimes_{i=1}^m \sO_{\e_i}.$$

The point of this construction is that twisting a line bundle by $\sO_{\e}$ 
changes it multidegree in the same manner as twisting by $(e,v)$, so that
if $\sL$ has multidegree $w$, then $\sL \otimes \sO_{w,w'}$ has multidegree
$w'$.

\begin{notn}\label{not:more-twist}
In Situation \ref{sit:geom}, suppose $\sL$ is a line bundle on
$\widetilde{X}$ of multidegree $w_0$. Then for any $w \in V(\bar{G}(w_0))$,
set
$$\sL_w := \sL \otimes \sO_{w_0,w}.$$

For $v \in V(\Gamma)$, set
$$\sL^v:=\sL_{w_v}|_{Z_v}.$$
\end{notn}

The following divisor sequences derived from our combinatorial data will
provide the backdrop to the multivanishing sequences and gluing conditions
considered in our definition of limit linear series.

\begin{notn}\label{not:twist-divs} In Situation \ref{sit:geom},
for each pair $(e,v)$ of an
edge and adjacent vertex of $\bar{\Gamma}$, let 
$D^{(e,v)}_0,\dots,D^{(e,v)}_{b_{v,v'}+1}$ be the sequence of effective 
divisors on $Z_v$ defined by $D^{(e,v)}_0=0$, and for $i\geq 0$,
$$D^{(e,v)}_{i+1}-D^{(e,v)}_{i}
=\sum_{\scriptsize \begin{matrix}\tilde{e}\text{ over }e:\\ 
\sigma(\tilde{e},v)\mu_v(\tilde{e}) \equiv 
- i \pmod{\bn(\tilde{e})}\end{matrix}} 
P_{\tilde{e}},$$
where $P_{\tilde{e}}$ denotes the node of $X$ corresponding to $\tilde{e}$.
\end{notn}

A global line bundle, together with our twisting bundles, induce gluing
isomorphisms as follows.

\begin{prop}\label{prop:gluing-maps} Given $\sL$ on $\widetilde{X}$ of
multidegree $w_0$, and vertices $v,v'$ of $\bar{\Gamma}$ connected by an
edge $e$, then for $i=0,\dots,b_{v,v'}$ we have isomorphisms
$$\vp^{(e,v)}_i:\sL^v(-D^{(e,v)}_i)/\sL^v(-D^{(e,v)}_{i+1}) \risom
\sL^{v'}(-D^{(e,v')}_{b_{v,v'}-i})/\sL^{v'}(-D^{(e,v')}_{b_{v,v'}+1-i})$$
induced by the line bundle $\sL_{w(v,v',i)}$, where 
$w(v,v',i) \in V(\bar{G}(w_0))$ is the $i$th vertex between $w_v$ and 
$w_{v'}$.
\end{prop}

This is essentially the last part of Proposition 4.4 of \cite{os25}.

\begin{defn}\label{def:critical} If $D_0,\dots,D_{b+1}$ is a non-decreasing
sequence of effective divisors on a smooth proper curve, we say $j$
is \textbf{critical} for $D_{\bullet}$ if $D_{j+1} \neq D_j$.
\end{defn}

We can now give our definition of limit linear series.

\begin{defn}\label{defn:lls} In Situation \ref{sit:geom}, suppose we have
a tuple $(\sL,(V^v)_{v \in V(\Gamma)})$ with $\sL$ a line bundle of 
multidegree $w_0$ on $\widetilde{X}$, and each $V^v$ an $(r+1)$-dimensional
space of global sections of the resulting $\sL^v$.
For each pair $(e,v)$ in $\bar{\Gamma}$, let $a^{(e,v)}_0,\dots,a^{(e,v)}_r$
be the multivanishing sequence of $V^v$ along $D^{(e,v)}_{\bullet}$.
Then $(\sL,(V^v)_{v \in V(\Gamma)})$ is a \textbf{limit linear series} if
for any $e \in E(\Gamma)$, with adjacent vertices $v,v'$, we have:
\begin{Ilist}
\itm for $\ell=0,\dots,r$, if $a^{(e,v)}_{\ell}=\deg D^{(e,v)}_j$ with $j$ 
critical for $D^{(e,v)}_{\bullet}$, then
\begin{equation}\label{eq:eh-genl}
a^{(e,v')}_{r-\ell} \geq \deg D^{(e,v')}_{b_{v,v'}-j};
\end{equation}
\itm there exist bases $s^{(e,v)}_0,\dots,s^{(e,v)}_r$ of $V^v$ and
$s^{(e,v')}_0,\dots,s^{(e,v')}_r$ of $V^{v'}$ such that
$$\ord_{D^{(e,v)}_{\bullet}}s^{(e,v)}_{\ell} = a^{(e,v)}_{\ell},
 \quad \text{ for } \ell=0,\dots,r,$$
and similarly for $s^{(e,v')}_{\ell}$,
and for all $\ell$ with \eqref{eq:eh-genl} an equality, we have
$$\vp^{(e,v)}_{j}(s^{(e,v)}_{\ell})=s^{(e,v')}_{r-\ell}$$
when we consider $s^{(e,v)}_{\ell} \in V^v(-D^{(e,v)}_{j})$ and 
$s^{(e,v')}_{r-\ell} \in V^{v'}(-D^{(e,v')}_{b_{v,v'}-j})$, where 
$j$ is as in (I), and $\vp^{(e,v)}_j$ is
as in Proposition \ref{prop:gluing-maps}.
\end{Ilist}

We say a limit linear series is \textbf{refined} if \eqref{eq:eh-genl} holds
with equality for all $\ell$.
\end{defn}

\begin{notn} In Situation \ref{sit:geom}, let
$$G^r_{\bar{w}_0}(X,\bn,(\sO_v)_v)$$
denote the moduli scheme of limit linear series.
\end{notn}

The notation $\bar{w}_0$ reflects that the space 
$G^r_{\bar{w}_0}(X,\bn,(\sO_v)_v)$ depends on $w_0$ only up to
arbitrary twists. For the construction of this moduli scheme, see
\S 3 of \cite{os25}. 

\begin{rem}\label{rem:concen-diff} Although our definition of concentrated
is more restrictive than that of \cite{os25}, this does not cause any
technical difficulties. Indeed, it is easy to see that tuples $(w_v)_v$
as in Situation \ref{sit:concen} always exist despite our more 
restrictive definition, and Proposition 3.5 of \cite{os25} asserts that the 
notion of limit linear series is in fact independent of the choice of tuple 
$(w_v)_v$. However, we briefly explain the relationship between our present 
definition and the definition in \cite{os25}. 

If we wanted an equivalent definition
of concentrated to that given in \cite{os25}, we could say that
$w$ is concentrated at a vertex $v \in V(\Gamma)$ if there is an ordering 
$$V(\Gamma)=\{v_1,v_2,\dots\}$$ 
with $v=v_1$, and such that for each $i>1$, we have that $w$
becomes negative in index $v_i$ after taking the composition of the 
twists at $(e_j,v_i)$ over all $j<i$ with $v_j$ adjacent to $v_i$, where 
$e_j$ is the edge of $\bar{\Gamma}$ connecting $v_i$ to $v_j$.
Indeed, this is equivalent to the definition in \cite{os25} because the 
degree in index $v_i$ is the same after the above-described twists as after 
taking the negative twists at all $v_j$ for $j<i$.
We thus see that the definition we are using is more restrictive, as claimed,
since it is equivalent to considering an ordering which is consistent with
the distance from $v$ in $\bar{\Gamma}$.
\end{rem}

\subsection*{Relationship to tropical Brill-Noether theory} 
Although our theory of limit linear series is quite different from the 
theory of divisors on graphs, certain definitions translate directly, as
we now explain.
 
First, our chain structures correspond to lengths of edges of metric 
graphs, with the distinction that we restrict to integer lengths. More
specifically, if we have a one-parameter smoothing of our nodal curve $X$,
then both the chain structure and the associated metric graph are 
determined by looking at the number of exceptional components lying over
each node of $X$ in a resolution to a regular smoothing.

Next, our notion of admissible multidegree corresponds to divisors on
the graph $\widetilde{\Gamma}$, with the added constraint that the divisor
must be nonnegative on all vertices of $\widetilde{\Gamma}$ lying over
edges of $\Gamma$, and the vertices over a given edge can have total degree
at most $1$.

Our notion of an admissible multidegree being concentrated at $v$ then
corresponds roughly to a $v$-reduced divisor on $\widetilde{\Gamma}$ --
a $v$-reduced divisor is concentrated at $v$, but not necessarily 
conversely (for instance, we do not require nonnegativity).

Finally, our twists of multidegrees at $(e,v)$ correspond to chip-firings 
along the set of vertices on the same side of $e$ as $v$ is. Thus, our
$\bar{G}(w_0)$ consists of a subset of the divisors on $\widetilde{\Gamma}$
linearly equivalent to $w_0$.

\section{Linear series with imposed multivanishing}\label{sec:mult-van}

In this section, we prove a generalized Brill-Noether theorem for
linear series with prescribed multivanishing, generalizing the theorem
of Eisenbud and Harris (Theorem 4.5 of \cite{e-h1}) from the case of
usual vanishing sequences to the case of multivanishing sequences.
We prove the genus-$0$ case with at most two multivanishing sequences
directly, and then use limit linear series for curves of compact type
to prove the main theorem.

We begin by setting up some basic notation and definitions.

\begin{notn}\label{not:imposed-multivan}
Let $X$ be a smooth projective curve of
genus $g$, and fix integers $r,d,n > 0$, and for $i=1,\dots,n$ fix also 
nondecreasing sequences $D^i_{\bullet}$ of effective divisors on $X$, such 
that the support of $D^i_{\bullet}$ is disjoint from that of
$D^{i'}_{\bullet}$ for every $i \neq i'$. Fix also a tuple of
nondecreasing sequences $a^i$, such that for each $i,j$ we have 
$a^i_j=\deg D^i_{\ell}$ for some $\ell$ critical for $D^i_{\bullet}$,
and the number of repetitions of $a^i_j$ is at most 
$\deg D^i_{\ell+1}- D^i_{\ell}$.

Then we denote by $G^r_d(X,(D^i_{\bullet},a^i)_i)$ the space of $\fg^r_d$s 
on $X$ having multivanishing sequence at least $a^i$ along $D^i_{\bullet}$ 
for each $i$.
\end{notn}

\begin{defn}\label{defn:strong-genl} Let $X$ be a smooth projective curve of
genus $g$, and fix integers $r,d,n > 0$, and for $i=1,\dots,n$ fix also 
$m_i>0$. Choose distinct points $P^i_j$ on $X$ for $i=1,\dots,n$ and 
$j=1,\dots,m_i$. Then we say that $(X,(P^i_j)_{i,j})$ is \textbf{strongly
Brill-Noether general} for $r,d$ if, for all tuples of nondecreasing 
effective divisor sequences $D^i_{\bullet}$, such that every divisor in 
$D^i_{\bullet}$ is supported among $P^i_1,\dots,P^i_{m_i}$,
and for every tuple of nondecreasing sequences $a^i$ as in Notation
\ref{not:imposed-multivan}, the space $G^r_d(X,(D^i_{\bullet},a^i)_i)$
has the expected dimension
$$\rho:=g+(r+1)(d-r-g)-\sum_{i=1}^n \left(\sum_{j=0}^r (a^i_j-j)+
\sum_{\ell=0}^{b_i} \binom{r^i_{\ell}}{2}\right)$$
if it is nonempty.

In the above, $D^i_{\bullet}$ is indexed from $0$ to $b_i+1$, and 
$r^i_{\ell}$ is defined to be $0$ if $\ell$ is not critical for 
$D^i_{\bullet}$, and the number of times $\deg D^i_{\ell}$ occurs in
$a^i$ if $\ell$ is critical.
\end{defn}

We will sometimes refer to an $X$ together with choices of $P^i_j$ as above 
as a `multimarked curve.'

The following is then the main theorem of this section.
We also take the opportunity to state some specific situations in which
Brill-Noether generality is known under precisely stated conditions.

\begin{thm}\label{thm:bn-genl} With notation as in Definition
\ref{defn:strong-genl}, we have that $(X,(P^i_j)_{i,j})$ is
strongly Brill-Noether general for $r$ and $d$ if any of the
following conditions are satisfied:
\begin{Ilist} 
\itm either $\ch k=0$ or $\ch k = p>d$, and both $X$ and the $P^i_j$
are general;
\itm $n \leq 2$ and both $X$ and the $P^i_j$ are general;
\itm $g=0$ and $n \leq 2$;
\itm either $\ch k=0$ or $\ch k = p>d$, $g=0$, and all $m_i$ are equal
to $1$;
\itm $g=1$, $n \leq 2$, $m_i=1$ for each $i$, and if $n=2$, then 
$P^1_1-P^1_2$ is not $\ell$-torsion for any $\ell \leq d$;
\itm $g=2$, $n \leq 1$, and if $n=1$, then $m_1=1$ and $P_1$ is not
a Weierstrass point.
\end{Ilist}
\end{thm}

Cases (I) and (II) should be considered the main result of the theorem.
Cases (IV)-(VI) fall into the case of classical vanishing
sequences, and thus were already known. The reason we restate them here
is that they will substantially broaden the number of explicit reducible
curves for which we can show that limit linear series spaces have the
expected dimension in Corollary \ref{cor:few-nodes-genl} below. 
We also mention that the $g=0$ case of (I) was
already proved by Garc\'ia-Puente \textit{et.\ al}.\ in \cite{g-h-h-m-r-s-t}.

Note also that our characteristic hypotheses are essentially optimal, as 
even in the special case treated by Eisenbud and Harris, the characteristic 
hypotheses are already necessary, even in genus $0$ for $r=1$.

\begin{proof} We first consider case (III).
Of course, if $n=0$ there is nothing to show, and if $n=1$ the imposition
of the multivanishing sequence $a^1$ is just a Schubert cycle in 
$G(r+1,\Gamma(X,\sO(d)))$, of codimension
$\sum_{j=0}^r (a^1_j-j)+ \sum_{\ell=0}^{b_1} \binom{r^1_{\ell}}{2}$.
Indeed, such a Schubert cycle for a partial flag is the same as the
Schubert cycle obtained from any completion of the flag, where any 
repetitions $a_{j_1}=a_{j_1+1}=\dots=a_{j_2}$ are replaced by
$a_{j_1},a_{j_1}+1,\dots,a_{j_1}+j_2-j_1$.
Similarly, in the case $n=2$, because the support of $D^1_{\bullet}$ is 
disjoint from that of $D^2_{\bullet}$, we see that the corresponding flags 
meet transversely, 
so the associated Schubert cycles intersect in the expected dimension,
as desired.

We now prove cases (I) and (II). If $X$ is general, then for $n$ 
general marked points $Q_1,\dots,Q_n$, under our hypotheses we have that 
the space of linear series on $X$ with prescribed ramification at the
$Q_i$ has the expected dimension; see for instance \cite{os18} (but note
that the proof of the main theorem can be simplified now that it is 
possible to construct a proper moduli space of limit linear series in
smoothing families, as carried out in \cite{os20}).
Accordingly, consider a one-parameter family in which for each $i$, all
the points in the support of $D^i_{\bullet}$ approach $Q_i$. Blowing up the 
special fiber at the $Q_i$, we obtain a curve $X_0$ consisting of a copy of 
$X$, together with rational tails glued at each $Q_i$, and with each 
$D^i_{\bullet}$
specializing onto the corresponding rational tail. But by construction
and the previously addressed rational case, we see that the space of 
limit linear series on $X_0$ with imposed multivanishing at the (limits of
the) $D^i_{\bullet}$ has the desired dimension, and the theorem follows.

As mentioned above, cases (IV)-(VI) were previously known: see for
instance Theorem 2.3 of \cite{e-h4}, Lemma 2.1 of \cite{os18}, and 
Theorem 1.1 of \cite{e-h6}. Thus, the theorem is proved.
\end{proof}

\section{Expected dimension for limit linear series}\label{sec:exp-dim}

In a sense, the smoothing theorem (Theorem 6.1)
of \cite{os25} already says that the ``expected dimension''
of the space of limit linear series on a curve of pseudocompact type
is $\rho$, in the sense that it is at least $\rho$, and if the conditions
cutting it out have maximal codimension, then it is exactly $\rho$.
However, these conditions -- ``linked determinantal loci'' -- are rather
abstract, so it is the goal of the present section to use the alternative
definition of limit linear series to give more geometrically concrete 
criteria for when the dimension is as expected. We then apply these
criteria in subsequent sections to more explicit families of curves. Our 
main result is as follows:

\begin{thm}\label{thm:exp-dim} In Situation \ref{sit:geom},
suppose further that each component of $X$
(considered as a multimarked curve) is strongly Brill-Noether general
for $r$ and $d$. Then the limit linear series space
$G^r_{\bar{w}_0}(X,\bn,(\sO_v)_v)$ has dimension $\rho$ if 
the gluing conditions imposed by Definition \ref{defn:lls} (II)
impose the maximal possible codimension. Furthermore, in this case
the refined limit linear series are dense.
\end{thm}

Recall that a number of sufficient conditions for strong Brill-Noether 
generality are listed in Theorem \ref{thm:bn-genl} above.

The main ingredient in Theorem
\ref{thm:exp-dim} is a combinatorial calculation which is the analogue
of the Eisenbud-Harris ``additivity of the Brill-Noether number'' 
(Proposition 4.6 of \cite{e-h1}). Despite the restriction to the 
pseudocompact-type case, the calculation is considerably more complicated, 
due to the presence of repeated vanishing orders and gluing conditions.

\begin{lem}\label{lem:combin-ineq} In the situation of Notation 
\ref{not:twist-divs}, let $Z_1$ and $Z_2$ be components of $X$ corresponding
to vertices $v,v'$ of $\bar{\Gamma}$ connected by an edge $e$. Write
$D^1_{\bullet}:=D^{(e,v)}_{\bullet}$ and 
$D^2_{\bullet}:=D^{(e,v')}_{\bullet}$,
let $\{0,\dots,b+1\}$ be the index set for $D^1_{\bullet}$, and let $C$ be 
the subset of critical indices. Observe that $\deg D^1_j+\deg D^2_{b+1-j}$
is independent of $j$, and denote its common value by $c$. For each $j$,
also set 
$$f_j:=\deg D^2_{c+1-j}-\deg D^2_{c-j}=\deg D^1_{j+1}-\deg D^1_j.$$
Given also sequences $a^1,a^2$ satisfying the conditions of 
Notation \ref{not:imposed-multivan}, let let $r^i_j$ be as in 
Definition \ref{defn:strong-genl}. Finally, set
$$g_j=r^1_j+\sum_{m \in C, m<j}\left(r^1_m-r^2_{b-m}\right).$$

We then have
\begin{equation}\label{eq:combin-ineq}
\sum_{j \in C} g_j(f_j+g_j-r^1_j-r^2_{b-j}) \geq 
(r+1)(c-1) - \sum_{\ell=0}^r (a^1_{\ell}+a^2_{r-\ell}) 
-\sum_{j \in C}\left(\binom{r^1_j}{2}+\binom{r^2_j}{2}\right),
\end{equation}
with equality precisely in the refined case.
\end{lem}

One checks using (4.4) from the proof of Lemma 4.6 of \cite{os25} that the 
$g_j$ defined above are precisely the number
of sections $s_{\ell}$ for which the gluing condition of Definition
\ref{defn:lls} (II) is imposed in index $j$.

\begin{proof} We first manipulate the righthand side of 
\eqref{eq:combin-ineq}. We have
\begin{align}\label{eq:combin-simple}
\nonumber (r+1)&(c-1) -  \sum_{\ell=0}^r (a^1_{\ell}+a^2_{r-\ell}) 
-\sum_{j \in C}\left(\binom{r^1_j}{2}+\binom{r^2_j}{2}\right) \\
\nonumber & = (r+1)(c-1) - \sum_{j \in C} \left(r^1_j \deg D^1_j 
+ r^2_{b-j} \deg D^2_{b-j}\right)
-\sum_{j \in C}\left(\binom{r^1_j}{2}+\binom{r^2_j}{2}\right) \\
& = (r+1)(c-1) - \sum_{j \in C} \left(r^1_j \left(c-f_j-\deg D^2_{b-j}\right)
+ r^2_{b-j} \deg D^2_{b-j}\right)\\
\nonumber & \quad\quad\quad \quad\quad\quad \quad
 +\sum_{j \in C} \left(\frac{\left(r^1_j\right)^2}{2}
-\frac{r^1_j}{2} +\frac{\left(r^2_{b-j}\right)^2}{2}
-\frac{r^2_{b-j}}{2}\right) \\
\nonumber & = \sum_{j \in C} \left(f_j r^1_j 
+\deg D^2_{b-j}\left(r^1_j-r^2_{b-j}\right)
-\frac{\left(r^1_j\right)^2}{2} -\frac{\left(r^2_{b-j}\right)^2}{2}\right).
\end{align}

First, in the refined case, we have $r^1_j=r^2_{b-j}=g_j$ for all $j \in C$,
so we see from \eqref{eq:combin-simple} that both sides 
of \eqref{eq:combin-ineq} are equal to
$$\sum_{j \in C} \left(r^1_j f_j-\left(r^1_j\right)^2\right),$$
and in fact \eqref{eq:combin-ineq} is an identity. We next reduce to the
refined case.

In the general case, we claim that if the sequences $a^1,a^2$ are
not refined, then we can always decrease one of the $a^2_{\ell}$ and still
have an allowable sequence. Indeed, if $\ell_0$ is
maximal such that \eqref{eq:eh-genl}
is strict, then suppose that $a^1_{\ell_0}=\deg D^1_{j_0}$,
and $a^2_{r-\ell_0}=\deg D^2_{b-j_1}$, with $j_1,j_0 \in C$, and $j_1<j_0$
by hypothesis. Let $j_2 \in C$ be minimal greater than $j_1$. Then our 
specific claim is
that if we set $\hat{a}^2_{r-\ell}$ to be equal to $a^2_{r-\ell}$ for 
$\ell \neq \ell_0$, and $\hat{a}^2_{r-\ell_0}=\deg D^2_{b-j_2}$, then
we still have a valid sequence. Certainly, \eqref{eq:eh-genl}
will still be satisfied, and $\hat{a}^2$ is
nondecreasing, so it is enough to see that $\hat{a}^2_{r-\ell_0}$ does
not contain too many repetitions of $\deg D^2_{b-j_2}$. Now, the maximal
number of allowed repetitions of $\deg D^2_{b-j_2}$ is 
$$\deg D^2_{b+1-j_2}-\deg D^2_{b-j_2}=\deg D^1_{j_2+1}-\deg D^1_{j_2}
=f_{j_2}.$$
Certainly, if $a^2_{r-\ell}=\deg D^2_{b-j_2}$, then $\ell>\ell_0$; by our
hypothesis for the maximality of $\ell_0$, for all such $\ell$ we have
$a^1_{\ell}= \deg D^1_{j_2}$. Now, if $j_0>j_2$, we see that in fact
there is no such $\ell$, because $a^1$ is nondecreasing. On the other hand,
if $j_0=j_2$, then we must have $a^1_{\ell_0+f_{j_0}}>\deg D^1_{j_0}$,
so by the maximality of $\ell_0$, we have 
$a^2_{r-\ell_0-f_{j_0}}<\deg D^2_{b-j_0}=\deg D^2_{b-j_2}$, so we have
at most $f_{j_0}-1=f_{j_2}-1$ repetitions of $\deg D^2_{b-j_2}$ in $a^2$, as 
desired.

Now, we examine the effect of replacing $a^2$ by $\hat{a}^2$ on both
sides of \eqref{eq:combin-ineq}. Clearly, the $r^1_j$ are unaffected,
while we are decreasing $r^2_{b-j_1}$ by $1$, and increasing $r^2_{b-j_2}$
by $1$. Recalling that $j_1$ and $j_2$ are consecutive
in $C$, we see that under our modification, $g_{j_2}$ increases by $1$, 
and the other $g_j$ are unchanged. From this, we compute that the
lefthand side of \eqref{eq:combin-ineq} increases by 
$g_{j_1}+g_{j_2}+f_{j_2}-r^1_{j_2}-r^2_{b-j_2}$, while using 
\eqref{eq:combin-simple}, the righthand side increases by 
$$\deg D^2_{b-j_1}-\deg D^2_{b-j_2}+r^2_{b-j_1}-r^2_{b-j_2}-1=
f_{j_2}+r^2_{b-j_1}-r^2_{b-j_2}-1.$$
Now, we have $g_{j_1} \leq r^2_{b-j_1}$ and $g_{j_2} \leq r^1_{j_2}$, but
in fact both inequalities are strict, by considering the corresponding
inequalities in our modified sequence.
We thus see that under our modification, the righthand side
of \eqref{eq:combin-ineq} increases by more than the lefthand
side. Now, after finitely many iterations, we will reach the refined
case, where we know both sides are equal, so we conclude that in general,
the righthand side is at most equal to the lefthand side, with equality
precisely in the refined case, as desired.
\end{proof}

\begin{proof}[Proof of Theorem \ref{thm:exp-dim}] By the definition of
strong Brill-Noether generality, it is enough to verify that
combinatorially, the expected dimension of 
$G^r_{w_{\bullet}}(X,\bn,(\sO_v)_v)$ is bounded by $\rho$, with strict
inequality for non-refined choices of multivanishing sequences. For 
$v \in V(\Gamma)$, let $d_v$ be the value of $w_v$ in index $v$, and
$g_v$ the genus of $Z_v$.
We stratify $G^r_{w_{\bullet}}(X,\bn,(\sO_v)_v)$ by multivanishing sequences,
so fix sequences $a^{(e,v)}$ for each adjacent pair $(e,v)$ in $\bar{\Gamma}$,
satisfying \eqref{eq:eh-genl}. Then 
$G^r_{w_{\bullet}}(X,\bn,(\sO_v)_v)$ has data consisting of a tuple of
$\fg^r_{d_v}$s on each $Z_v$ and gluing data for the underlying line bundle,
and imposed conditions are given by \eqref{eq:eh-genl} together with the 
gluing condition described in Definition \ref{defn:lls} (II).

The expected dimension of the data is
\begin{align} 
\nonumber \sum_v (g_v+(r+1)(d_v-r-g_v)) & + n
 =g+(r+1)\left(\sum_v d_v - r|V(\Gamma)|-g+n\right) \\
& =g+(r+1)\left(\sum_v d_v +|E(\Gamma)|- (r+1)|V(\Gamma)|-g+1\right) 
\end{align}
where $n=|E(\Gamma)|-|V(\Gamma)|+1$ is the number of parameters obtained
in gluing the line bundles on the $Z_v$ to obtain a line bundle on $X$, and
hence $n+\sum_v g_v = g$.

Now, we need to relate the $d_v$'s to $d$. Fix a $v$. By definition, 
$w_v$ has entry $d_v$ in index $v$. We also see that
for any $v' \neq v$, if $v''$ is the vertex adjacent to $v'$ in the
direction of $v$, and $n_{v',v''}$ denotes the number of edges of
$\Gamma$ connecting $v'$ to $v''$, then the entry of $w_v$ in index $v'$
plus the number of edges $e$ of $\Gamma$ connecting $v'$ to $v''$ such
that $\mu_v(e) \neq 0$ is equal to
$$d_{v'}-c_{v',v''}+n_{v',v''},$$
where $c_{v',v''}$ is obtained as the $c$ of Lemma \ref{lem:combin-ineq}
when we consider $Z_{v'}$ and $Z_{v''}$.
By definition, if we sum over all $v'$, we get $d$, so we conclude that
\begin{equation}\label{eq:deg} 
d = \sum_{v \in V(\Gamma)} d_v - \sum_{e \in E(\bar{\Gamma})} c_{h(e),t(e)}
+ |E(\Gamma)|.
\end{equation}
Thus, we can re-express our earlier expected dimension for the limit linear
series data as
\begin{multline*}
g+(r+1)\left(d +\sum_{e \in E(\bar{\Gamma})} c_{h(e),t(e)}
- (r+1)|V(\Gamma)|-g+1\right) \\
=\rho+(r+1)\sum_{e \in E(\bar{\Gamma})} (c_{h(e),t(e)}-(r+1)).
\end{multline*}

We next consider the codimension imposed by the multivanishing sequences.
For each pair $v,v'$ of adjacent vertices connected by an edge $e$ of
$\bar{\Gamma}$, the imposed multivanishing on $Z_v$ and $Z_{v'}$ at the
nodes corresponding to $e$ has codimension
\begin{multline*}
\sum_{\ell=0}^r (a^{(e,v)}_{\ell}-\ell+a^{(e,v')}_{r-\ell} -(r-\ell))
+\sum_{j=0}^{b_{v,v'}} \left(\binom{r^{(e,v)}_j}{2} 
+ \binom{r^{(e,v')}_j}{2}\right) \\
= -r(r+1) + \sum_{\ell=0}^r (a^{(e,v)}_{\ell}+a^{(e,v')}_{r-\ell}) 
+\sum_{j=0}^{b_{v,v'}} \left(\binom{r^{(e,v)}_j}{2} 
+ \binom{r^{(e,v')}_j}{2}\right).
\end{multline*}
Subtracting this from our previous expected dimension, we find that
what remains is
$$\rho+\sum_{e \in E(\bar{\Gamma})} \left((r+1)(c_{h(e),t(e)}-1)
- \sum_{\ell=0}^r (a^{(e,v)}_{\ell}+a^{(e,v')}_{r-\ell}) 
-\sum_{j=0}^{b_{v,v'}} \left(\binom{r^{(e,v)}_j}{2} 
+ \binom{r^{(e,v')}_j}{2}\right)\right).$$
But note that this is greater than $\rho$ by precisely the sum over
$e \in E(\bar{\Gamma})$ of the righthand side of \eqref{eq:combin-ineq}.
In order to confirm that we have the correct expected dimension, it
therefore suffices to confirm that the expected codimension of the gluing 
conditions is at least the lefthand side of \eqref{eq:combin-ineq}. In 
fact, we show that they are the same, with the exception that any negative 
summands in the lefthand side must be replaced by $0$.

But this is clear: each gluing condition is of the form that two subspaces,
of dimensions $r^1_j$ and $r^2_{b-j}$ respectively, of an $f_j$-dimensional
vector space, must intersect in dimension at least $g_j$. That is, it is
the preimage of a closed subset 
$Z \subseteq G(r^1_j,f_j) \times G(r^2_{b-j},f_j)$ which one easily verifies
is irreducible of codimension $g_j(f_j+g_j-r^1_j-r^2_{b-j})$ when the latter
is nonnegative, 
yielding the desired statement.
\end{proof}

\section{Curves with few nodes: restricted chain 
structures}\label{sec:few-nodes-1}

We now consider curves which are close to being of compact
type, in the sense that any two components intersect in at most three
nodes. For such curves, individual gluing conditions are straightforward to 
understand, so in light of Theorem \ref{thm:exp-dim}, if we also impose
strong Brill-Noether generality on the individual components, in order to 
prove the entire curve is Brill-Noether general, it is enough to show that 
the gluing conditions are independent of one another.
In this section, we show that suitable conditions on the chain structures
always imply the desired independence, irrespective of enriched structures.
The families we consider generalize the curves studied by Cool, Draisma, 
Payne and Robeva in \cite{c-d-p-r}, and our results may be viewed as a 
generalized analogue of theirs, which simultaneously lends a geometric
interpretation to their purely numerical `genericity' condition. See
Remark \ref{rem:cdpr} for details.

In contrast, in the next section we will produce smaller families which we 
can show to be Brill-Noether general for arbitrary chain structures, under 
a generality hypothesis instead on the enriched structures.

We begin with a background proposition which holds for
arbitrary curves of pseudocompact type.

\begin{prop}\label{prop:concen-nonneg} The $(w_v)_v$ of Situation
\ref{sit:concen} can be chosen so that each $w_v$ is nonnegative in every
index. 
\end{prop}

\begin{proof} Let $v_1$ be any vertex of $\bar{\Gamma}$, and write 
$V(\bar{\Gamma})=\{v_1,v_2,\dots,v_m\}$ with the
ordering compatible with the distance from $v_1$ in $\bar{\Gamma}$.
If $w'_{v_1}$ is any admissible multidegree concentrated at $v_1$,
we can modify it into a $w_{v_1}$ which is nonnegative in all indices
other than $v_1$ as follows: for each $i>1$, let $e_i \in E(\bar{\Gamma})$
be the edge adjacent to $v_i$ in the direction of $v_1$, and $v'_i$ the
other vertex adjacent to $e_i$. If $w'_{v_1}$ is not nonnegative in
index $v_m$, twist $w'_{v_1}$ at $(e_m,v'_m)$ the minimal number of
times to make it nonnegative in index $v_m$. The minimality implies that
the new $w'_{v_1}$ is still concentrated at $v_1$. Then repeat the process
with $v_{m-1}$, $v_{m-2}$, and so forth, and ultimately we will arrive at
the desired $w_{v_1}$. 
We then obtain $w_{v_2}$ from $w_{v_1}$ by twisting at $(e_2,v_1)$ the
maximal number of times possible without making it negative in index $v_1$.
We continue in this manner, inductively obtaining $w_{v_i}$ from $w_{v'_i}$ 
by twisting at $(e_i,v'_i)$ the maximal number of times possible without
making it negative in index $v'_i$ (note that we always have $v'_i=v_j$ for
some $j<i$).
\end{proof}

Our expected dimension result for restricted chain structures is then
the following.

\begin{cor}\label{cor:few-nodes} In Situation \ref{sit:geom},
suppose further that
\begin{Ilist}
\itm there are at most three edges of $\Gamma$ connecting any given pair of
vertices;
\itm for any adjacent vertices $v,v'$ of $\Gamma$, if $v,v'$ are connected
by edges $(e_i)_i$, then for any integers $(x_i)_i$ with 
$\sum_i x_i \bn(e_i)=0$, if there is a unique $j$ with $x_j>0$, then 
we have 
$$\sum_i \left\lfloor x_j \bn(e_j)/\bn(e_i)\right\rfloor > d;$$ 
\itm each (multi)marked component of $X$ is strongly Brill-Noether general.
\end{Ilist}

Then the space of limit linear series on $(X,\bn)$ of degree $d$ is pure of 
the expected dimension $\rho$.
\end{cor}

\begin{rem}\label{rem:few-nodes-conds}
Note that condition (II) above is vacuous when $v,v'$ are connected by a 
single edge. When $v,v'$ are connected by a pair of edges $e_1,e_2$, then
condition (II) amounts to requiring that
$$\frac{\lcm(\bn(e_1),\bn(e_2))}{\bn(e_1)} +
\frac{\lcm(\bn(e_1),\bn(e_2))}{\bn(e_2)} > d.$$
See Remark \ref{rem:cdpr} below for the relationship to \cite{c-d-p-r}.

Note also that even when $v,v'$ are connected by three edges, it is
easy to find values of the $\bn(e_i)$ for which (II) is satisfied: for 
instance, one can take $1,n,n^2$ for any $n\geq d$,
or $1,n,n'$ for $n,n'$ relatively prime and at least equal to $d$.

Recall also that Theorem \ref{thm:bn-genl} (III)-(VI) gives explicit cases
in which condition (III) above is known to be satisfied.
\end{rem}

\begin{proof} According to Theorem \ref{thm:exp-dim}, condition (III) 
implies that it is enough to show
that the gluing conditions impose the maximum codimension. We begin by
describing what happens for a pair of nodes connected by two edges. 
Condition (I) implies that each nontrivial gluing condition is expected to 
impose codimension $1$; furthermore, such a condition can be imposed only when
$$\deg D^{(e,v)}_{i+1}-D^{(e,v)}_i=2$$
for some $i$, where $v$ is as in (II) and $e$ is the edge connecting $v$ to 
$v'$ in $\bar{\Gamma}$. Now, we claim that we
get at most one such nontrivial condition for any pair of components. 
Indeed, this is a consequence of condition (II) and Proposition 
\ref{prop:concen-nonneg}, since the way the $D^{(e,v)}_{\bullet}$
is defined, in order to have
$$\deg D^{(e,v)}_{i+1}-D^{(e,v)}_i=2,$$
we need to have
$\bn(e_1) | (\sigma(e_1,v) \mu_v(e_1)+i)$ and 
$\bn(e_2)|(\sigma(e_2,v)\mu_v(e_2)+i)$;
if this occurs 
for $i_1>i_2$, then $\lcm(\bn(e_1),\bn(e_2))|(i_1-i_2)$. But in order
for a given $\deg D^{(e,v)}_i$ to occur in the multivanishing sequence,
we must have $\deg D^{(e,v)}_i \leq d_v \leq d$,
where $d_v$ is the value of $w_v$ in index $v$. Now, 
\begin{align*} \deg D^{(e,v)}_{i_1}-\deg D^{(e,v)}_{i_2} 
& = \frac{i_1-i_2}{\bn(e_1)}+\frac{i_1-i_2}{\bn(e_2)} \\
& \geq \lcm(\bn(e_1),\bn(e_2))\left(1/\bn(e_1)+1/\bn(e_2)\right), 
\end{align*}
so (III) implies that this cannot occur, proving the claim.

Now, suppose $v$ and $v'$ are connected by edges $e_1$ and $e_2$, and we
are given linear series with the appropriate multivanishing on $Z_v$ and 
$Z_{v'}$. By the claim above, we have at most a single gluing condition at
$Z_v \cap Z_{v'}$, and there are two cases to consider: if the sections
determining the directions of gluing on $Z_v$ and $Z_{v'}$ are nonvanishing 
at both $P_{e_1}$ and $P_{e_2}$, then the gluing condition imposes the
desired codimension $1$ on the choice of global line bundle gluing at
$Z_v \cap Z_{v'}$. Otherwise, the given gluing condition can be considered
degenerate, since it is either impossible to glue, or every
gluing of line bundles satisfies the required gluing condition. We
thus stratify the space of tuples of linear series on the $Z_v$ having
the required multivanishing sequences based on which gluing conditions are
degenerate. The nondegenerate cases all
visibly independently impose codimension $1$ on the choice of gluing for
the global line bundle, while we see that for every degenerate gluing
condition our stratum itself has its dimension reduced by (at least) 
$1$,
again independently. Indeed, the condition of having an extra order
of vanishing at $P_{e_1}$ or $P_{e_2}$ can be expressed with a modified
multivanishing sequence by the insertion of an addition divisor; the
expected codimension of the modified multivanishing sequence is $1$ greater
than before, and by our strong Brill-Noether generality hypothesis, this
additional codimension is in fact realized, as desired.

Now we consider what happens when some $v,v'$ is connected by three edges
$e_1,e_2,e_3$. The argument is similar to the above,
except that the analysis of gluing conditions is
somewhat more complicated. Let $e$ denote the edge connecting
$v$ and $v'$ in $\bar{\Gamma}$; a nontrivial gluing condition may 
impose codimension $1$ or codimension $2$, depending partly on whether
$$D^{(e,v)}_{i+1}-D^{(e,v)}_i$$
has degree $2$ or $3$.
Now, we claim that the
inequalities of condition (II) imply that the latter case can occur at
most once, and the former case can occur at most twice. Moreover, if the 
latter case occurs, then the former case does not, and if the former case
does occur twice, the corresponding (reduced) degree-$2$ divisors cannot be 
supported on the same pair of points. Indeed, most of these assertions 
follow from the case that one of the $x_j$ is zero by the same argument 
as above, which implies that there cannot be any pair of
indices $i$ such that $D^{(e,v)}_{i+1}-D^{(e,v)}_i$
contains a given pair of points $P_{e_j}$ and
$P_{e_{j'}}$. 
The only difference is that with three nodes, the relevant formula for the
change in divisor degrees is
$$\deg D^{(e,v)}_{i_1}-\deg D^{(e,v)}_{i_2} 
\geq \sum_{j=1}^3 \left\lfloor \frac{i_1-i_2}{\bn(e_j)}\right\rfloor.$$
The only remaining possibility to consider is that 
$D^{(e,v)}_{i+1}-D^{(e,v)}_i$ has degree $2$ three times, supported at all
three possible pairs of points. But this is ruled out by condition (II)
in the case that all $x_j$ are nonzero.

Now, we proceed as before. Degeneracy can be more complicated if we have
$$\deg D^{(e,v)}_{i+1}-D^{(e,v)}_i=3:$$
if the relevant multivanishing
index occurs once, the corresponding one-dimensional space might vanish
at one or two of the points in $D^{(e,v)}_{i+1}-D^{(e,v)}_i$; if the 
relevant multivanishing index occurs
twice, the corresponding two-dimensional space might vanish at one of the
points, or it might contain a one-dimensional subspace vanishing at two
of the points. All of these cases can be treated via the strong
Brill-Noether generality hypothesis as before, and under hypothesis (IV) we 
see that they each impose the expected codimension, which is respectively 
$1$, $2$, $2$ and $1$.
Gluing conditions can also be more complicated: depending on the 
number of repetitions in complementary terms of the multivanishing 
sequences on $Z_v$ and $Z_{v'}$, 
we could obtain nontrivial gluing conditions in terms of agreement of 
a pair of points or a pair of lines in the projective
plane, or in terms of a point lying on a line; the first two conditions
have expected codimension $2$, while the latter has expected codimension 
$1$. However, in all combinations of degeneracy type and gluing conditions,
we see that either a given gluing becomes impossible, or the degeneracy
on both $Z_v$ and $Z_{v'}$ imposes enough codimension to make up for lost 
gluing conditions.
Thus, we obtain the desired gluing codimension in the case that 
$$\deg D^{(e,v)}_{i+1}-D^{(e,v)}_i=3$$
for some $i$.

The only other new case is that $D^{(e,v)}_{i+1}-D^{(e,v)}_i$
has degree $2$ twice for a given $v,v'$.
Here again the expected gluing codimension
is $2$, and if there is any degeneracy, the same analysis as before
gives the desired codimension. On the other hand, in the nondegenerate
case, one easily sees that the pair of nontrivial gluing conditions 
(supported on two distinct pair of points) uniquely determines the gluing
of the line bundles, and hence imposes codimension $2$, as desired.
\end{proof}

\begin{rem}\label{rem:cdpr}
If we set $d=2g-2$ in condition (II) of Corollary \ref{cor:few-nodes},
then in light of Remark \ref{rem:few-nodes-conds} and the discussion at 
the end of 
\S \ref{sec:review}, this recovers precisely the `genericity' condition of
\cite{c-d-p-r}. Thus, our proof of Corollary \ref{cor:few-nodes} is in
essence showing that for the curves they consider, and under their
genericity conditions, gluing conditions always automatically impose 
the expected codimension, and do so purely on the level of line bundle
gluings. It seems reasonable to expect that this sort of situation should 
be highly predictive for when metric graphs will be Brill-Noether general.
\end{rem}

We conclude this section with an explicit criterion for Brill-Noether
generality in terms of degenerations, simultaneously generalizing previous
conditions due to Eisenbud and Harris, Welters \cite{we4}, and Cools, 
Draisma, Payne and 
Robeva. For simplicity, we restrict to characteristic $0$, although this
only comes into condition (II), and the positive-characteristic analogue 
follows identically from Theorem \ref{thm:bn-genl}.

\begin{cor}\label{cor:few-nodes-genl} Let $\pi:X \to B$ be a flat, proper
morphism with $B$ the spectrum of a discrete valuation ring over a
field of characteristic $0$, generic
fiber $X_{\eta}$ a smooth curve, and special fiber $X_0$ a nodal curve.

Suppose that each component of $X_0$ has one of the following forms:
\begin{Ilist}
\itm a rational curve meeting at most two other components of $X_0$ in
at most three nodes each;
\itm a rational curve meeting every other component of $X_0$ in at most
one node;
\itm an elliptic curve meeting at most one other component of $X_0$ in
a single node, or meeting two other components of $X_0$ in single nodes
which do not differ by $m$-torsion for any $m \leq 2g-2$;
\itm a genus-$2$ curve meeting one other component of $X_0$ at a single
point which is not a Weierstrass point.
\end{Ilist}

Suppose further that the chain structure $\bn$ induced on $X_0$ by the 
singularities of $X$ satisfies the inequalities of condition (II) of
Corollary \ref{cor:few-nodes}, with $d=2g-2$.

Then $X_{\eta}$ is Brill-Noether general.
\end{cor}

\begin{proof} According to Theorem \ref{thm:bn-genl} (III)-(VI), the
hypotheses on the components of $X_0$ imply that they are strongly
Brill-Noether general, so condition (III) of Corollary \ref{cor:few-nodes}
is satisfied. Condition (I) is visibly satisfied, and we finally observe
that condition (II) is invariant under scaling of the chain structure,
so Corollary \ref{cor:few-nodes} tells us that we have the expected dimension
for all limit linear series spaces for $(X_0,c\bn)$, with $c$ any positive 
integer. It then follows from the specialization result Corollary 3.14
of \cite{os25} that $X_{\eta}$ is Brill-Noether general. 
\end{proof}

\section{Curves with few nodes: general enriched 
structures}\label{sec:few-nodes-2}

Continuing with the study of curves with few nodes between given pairs
of components, we consider a narrower family of curves, 
and show that every curve in the family is Brill-Noether general if 
equipped with a general enriched structure, irrespective of chain structure. 

We begin with a straightforward proposition on the structure of choices of 
enriched structures, which holds for arbitrary curves of pseudocompact type.

\begin{prop}\label{prop:enriched-structure} In Situation \ref{sit:geom},
the space of enriched structures on $(X,\bn)$ is canonically identified
with the product of the spaces of enriched structures on subcurves
$Z_{v,v'} \subseteq \widetilde{X}$, where $v,v'$ are adjacent vertices
of $\Gamma$, and $Z_{v,v'}$ is the subcurve of $\widetilde{X}$ containing
$Z_v$ and $Z_{v'}$, and all chains of rational components connecting 
them.

If $v,v' \in V(\Gamma)$ are connected by $m$ edges $e_1,\dots,e_m$, then
$(k^*)^{m-1}$ acts freely on the space of enriched structures on
$Z_{v,v'}$ as follows: given 
$(\lambda_1,\dots,\lambda_{m-1}) \in (k^*)^{m-1}$, we scale 
the gluing map defining $\sO_v$ at $P_{e_i}$ by $\lambda_i$ 
and the gluing map defining $\sO_{v'}$ at $P_{e_i}$ by $\lambda_i^{-1}$
for each $i \leq m-1$.
\end{prop}

Our expected dimension result is then the following.

\begin{cor}\label{cor:two-node-two} In Situation \ref{sit:geom},
suppose further that
\begin{Ilist}
\itm there are at most two edges of $\Gamma$ connecting any given pair of
vertices;
\itm for every $v \in V(\Gamma)$, if there is some $v' \in V(\Gamma)$
which is connected to $v$ by two edges, then $Z_v$ is rational and $v$
has valence at most three in $\Gamma$;
\itm the enriched structure on $(X,\bn)$ is general;
\itm each (multi)marked component of $X$ is strongly Brill-Noether general.
\end{Ilist}

Then the space of limit linear series on $(X,\bn)$ of degree $d$ is pure of 
the expected dimension $\rho$.
\end{cor}

Thus, the curves considered in Corollary \ref{cor:two-node-two} are 
almost of compact type, except that an elliptic component containing at
most two nodes may be replaced by a pair of rational components meeting
each other at a pair of nodes. Nonetheless, this contains (noncompact-type)
curves of every genus, some of which have arisen in other contexts; see 
Example \ref{ex:chain-hyperelliptic}.

The following lemma describes the behavior of gluing conditions in the
case of interest, saying in essence that with fixed (multi)vanishing
sequences, we will always have some gluing directions fixed, and the
rest varying freely.

\begin{lem}\label{lem:rational-gluing} Let $X$ be rational, with
distinct points $P_1,P_2$ and $Q$. Fix also divisor sequences
$D^1_{\bullet}, D^2_{\bullet}$ as in Definition \ref{def:multivanishing} 
and supported on the $P_i$ and $Q$ respectively, $r,d>0$, and multivanishing 
sequences 
$a^1$ and $a^2$, and let $G^{r,\circ}_d(X,(D^i_{\bullet},a^i)_i)$ denote the
open subscheme of $G^r_d(X,(D^i_{\bullet},a^i)_i)$ consisting of linear series
with multivanishing sequences along the $D^i_{\bullet}$ given precisely by 
$a^i$.  Let $S$ be the set of indices $j$ such
that $\deg D^1_{j+1}-\deg D^1_{j}=2$, and
$\deg D^1_j$ occurs exactly once in $a^1$. Then there exists a subset
$S' \subseteq S$ such that the image of the natural map 
\begin{equation}\label{eq:special-dirs}
G^{r,\circ}_d(X,(D^i_{\bullet},a^i)_i) \to 
\prod_{j \in S} \PP(\sO_X(d)(-D^1_j)/\sO_X(d)(-D^1_{j+1}))
\end{equation}
is a fiber of the projection to 
$$\prod_{j \in S'} \PP(\sO_X(d)(-D^1_j)/\sO_X(d)(-D^1_{j+1})).$$

Moreover, the image in
$$\prod_{j \in S'} \PP(\sO_X(d)(-D^1_j)/\sO_X(d)(-D^1_{j+1}))$$
is nondegenerate in the sense that it is a tuple of subspaces, each
of which has no more vanishing at either $P_1$ or $P_2$ than that
prescribed by $a^1$.

Finally, the induced map
\begin{equation}\label{eq:special-dirs-2}
G^{r,\circ}_d(X,(D^i_{\bullet},a^i)_i) \to 
\prod_{j \in S \smallsetminus S'} \PP(\sO_X(d)(-D^1_j)/\sO_X(d)(-D^1_{j+1}))
\end{equation}
has every fiber of codimension equal to $\#(S\smallsetminus S')$.
\end{lem}

\begin{proof} First observe that the space $G^{r}_d(X,(D^i_{\bullet},a^i)_i)$,
as the intersection of two Schubert varieties corresponding to 
complementary flags, is a Richardson variety and hence irreducible.
We refer to \S\S 2.2 and 2.5 (particularly Lemma 2.6) of Vakil \cite{va7} 
for background on
intersections of pairs of Schubert varieties, expressed combinatorially
in terms of positions of black and white checkers. By irreducibility, it is 
enough to prove the first statement of the lemma for any open subset of
$G^{r,\circ}_d(X,(D^i_{\bullet},a^i)_i)$. Let 
$G^{r,\circ\circ}_d(X,(D^i_{\bullet},a^i)_i)$ be the largest cell, 
consisting of 
spaces $V$ admitting a basis $e_0,\dots,e_r$ such that $e_i$ vanishes to 
order exactly $a^1_i$ along $D^1_{\bullet}$, and to order exactly 
$a^2_{r-i}$ along $D^2_{\bullet}$.
Now, in order for
$G^{r}_d(X,(D^i_{\bullet},a^i)_i)$ to be nonempty we must have
$a^1_i+a^2_{r-i} \leq d$ for all $i$; let $S'$ be the set of $j \in S$
such that we have $\deg D^1_j=a^1_i$ with
$a^1_i+a^2_{r-i} = d$, and let $T$ be the set of $i$ such that 
$a^1_i=\deg D^1_j$ for some $j \in S'$.

If $i \in T$, the choice of $e_i$ is unique up to scalar, and since $a^1_i$
is nonrepeated, the image of \eqref{eq:special-dirs} in index $a^1_i$ is
thus unique on $G^{r,\circ\circ}_d(X,(D^i_{\bullet},a^i)_i)$. Moreover, 
$e_i$ cannot
have any additional vanishing at $P_1$ or $P_2$. On the other hand,
if $i \not\in T$, then $h^0(X,\sO_X(d)(-D^1_j)) \geq 2$, and therefore the
choices of $e_i$ surject onto the nonzero elements of 
$\sO_X(d)(-D^1_j)/\sO_X(d)(-D^1_{j+1})$. Because the choices of $e_i$
are independent for different $i$, we conclude that the image of
$G^{r,\circ\circ}_d(X,(D^i_{\bullet},a^i)_i)$ under \eqref{eq:special-dirs} 
is a fiber of the projection to
$\prod_{j \in S'} \PP(\sO_X(d)(-D^1_j)/\sO_X(d)(-D^1_{j+1}))$,
as desired, and has the claimed nondegeneracy property.

The proof of the assertion on the fibers \eqref{eq:special-dirs-2} is 
similar, but a bit more involved: if we fix a point 
$$z \in \prod_{j \not \in S'} \PP(\sO_X(d)(-D^1_j)/\sO_X(d)(-D^1_{j+1})),$$
then the fiber over $z$ still corresponds to (an open subset of) an
intersection of two Schubert varieties, with one still corresponding to the
flag determined by $D^2_{\bullet}$, and the other corresponding to a 
refinement of the flag determined by $D^1_{\bullet}$, with exactly 
$\#\{S\smallsetminus S'\}$ (non-repeated) entries of $a^1$ increased by $1$. 
Let $(a^1)'$ denote this modified multivanishing sequence. The new flags 
need not be transverse, and the intersection may no longer be irreducible,
but we can study its maximal cells as described in \cite{va7}, and it
suffices to see that their dimensions are still as expected (i.e., are the 
same as in the case that the flags are transverse). Now, in Vakil's notation
the non-transversality of the refined flag is expressed in terms of black
checkers, by starting with the default position along the antidiagonal
(corresponding to two transverse flags), and allowing some (disjoint) pairs 
of adjacent checkers to have their $x$ and $y$-coordinates swapped.
Possible white checkers positions correspond to permutations $\sigma$ of
$\{0,\dots,r\}$, with the $(i+1)$st white checker in position 
$((a^1)'_{\sigma(i)}+1,a^2_{r-i}+1)$.
We also have the additional constraint
that each white checker must have a black checker weakly above it, and a
black checker weakly to the left of it. We first consider the case that
$\sigma=\id$, which corresponds to the maximal cell in the case of transverse
flags. We see that because $z$ only specified directions away from $S'$,
the entries of $(a^1)'$ which were increased by $1$ did not correspond to
white checkers which were already on the antidiagonal, so the new 
configuration of white checkers lies entirely on or below the antidiagonal.
In this case, we see that our adjacent swaps of black checkers does not 
affect the number of black checkers dominated by any given white one, 
so the dimension of this cell is as expected. Finally, given any $\sigma$,
if we write it as a product of disjoint cycles we see that the dimension
of the corresponding cell can be analyzed cycle by cycle, and any increase
in black checkers dominating by the white ones is always at least cancelled
by a corresponding increase in domination of other white checkers, so the
dimension cannot increase.
We thus conclude that we have the expected codimension, as desired.
\end{proof}

\begin{proof}[Proof of Corollary \ref{cor:two-node-two}]
According to Theorem \ref{thm:exp-dim}, it is enough to verify that
if we stratify the space of limit linear series by multivanishing 
sequences, the gluing conditions impose the maximal possible 
codimension. Note that because there are finitely many strata, it is 
enough to prove that a general enriched structure has the desired behavior
one stratum at a time. According to Proposition 
\ref{prop:enriched-structure}, we can describe enriched structures in terms
of one pair of adjacent components at a time. Now, if $v$ and $v'$ are 
connected by $2$ edges, then modifying the enriched structure by scaling
as in the proposition will change the resulting gluing maps by successive
powers of the given scalar. Thus, for any two fixed gluing directions on
$Z_v$ and on $Z_{v'}$, a general enriched structure will not allow both
conditions to be simultaneously satisfied. According to Lemma 
\ref{lem:rational-gluing}, with the multivanishing sequences fixed on
both $Z_v$ and $Z_{v'}$, we find a finite set of fixed gluing directions
on each of $Z_v$ and $Z_{v'}$, and the generality condition on the 
enriched structure simply requires that if two or more such fixed directions
on $Z_v$ are paired with those on $Z_{v'}$, that there not exist any
line bundle gluing simultaneously satisfying the gluing conditions. 

With this generality condition, we have ensured that if a given stratum
is nonempty, each pair $Z_v,Z_{v'}$ as above has at most one gluing
condition which comes from fixed directions on both sides. Such a 
condition, if it occurs, uniquely determines the line bundle gluing
between $Z_v$ and $Z_{v'}$, while any gluing condition which is not
fixed on one side or the other imposes the desired codimension, 
independently of one another, by the last part of Lemma 
\ref{lem:rational-gluing}. We thus conclude that the gluing conditions
impose the maximal codimension, as desired.
\end{proof}

\begin{ex}\label{ex:chain-hyperelliptic} The prototypical example of
the curves treated by Corollary \ref{cor:two-node-two} is a chain
of rational curves in which the components alternate meeting in two nodes
or one node. Such curves were considered for instance in Jensen and Payne's
tropical approach to the Gieseker-Petri theorem \cite{j-p1}.
Corollary \ref{cor:few-nodes} and 
\ref{cor:two-node-two} imply that such curves are Brill-Noether general
with suitable hypotheses on either the chain or enriched structure. 
However, we see that if arbitrary chain and enriched structures are allowed,
they are not Brill-Noether general. Indeed, such curves are in the closure 
of the hyperelliptic locus (for instance, one easily checks that they carry
degree-$2$ admissible covers of a chain of projective lines),
so general theory tells us that there must exist limit
$\fg^1_2$s on them. 

However, it is just as easy to see this explicitly: take the
trivial chain structure (meaning $\bn(e)=1$ for all $e$), and choose 
concentrated multidegrees of $2$
on the main component and $0$ on every other component. Take 
(multi)vanishing sequence $0,2$ at each node or pair of nodes. On each
interior component, we are looking for $2$-dimensional spaces of 
polynomials of degree $2$ which have an element vanishing to order $2$ at
one point, and another element vanishing simultaneously at two other
points. Such spaces are of course uniquely determined. At the ends, we
have only the condition that our space contain an element vanishing
simultaneously at two points, so have a projective line of choices. 
Now, there are no gluing conditions at the pairs of components connected
by a single node, but at the pairs connected by two nodes, we have 
two gluing conditions, with one coming from each term of the multivanishing 
sequence. The first gluing condition uniquely determined a choice of gluing 
for the line bundle, but the second gluing condition will only be satisfied 
for special enriched structures, except on the first and last components.
On the interior components, one checks easily that there does always exist 
enriched structures for which the two gluing conditions are compatible,
so we do in fact get a limit $\fg^1_2$, as claimed. 
\end{ex}

\section{Binary curves}\label{sec:binary}

As the opposite extreme from the case where there are few nodes connecting
any pair of components, we now consider the case of binary curves:

\begin{sit}\label{sit:binary} Suppose that $X$ is obtained by gluing
rational curves $Z_1$ and $Z_2$ to one another at $g+1$ nodes.
\end{sit}

In Situation \ref{sit:binary}, a $\fg^r_{(d_1,d_2)}$ refers to 
a line bundle on $X$ of degree $d_i$ on each $Z_i$, together with an
$(r+1)$-dimensional space of global sections.
In the case of binary curves, in contrast to the previous sections, it turns 
out that there is only a single gluing condition which can fail to impose 
maximal codimension, and this condition is insensitive to chain structures 
or enriched structures. We will show that for any binary curve, many 
components of the spaces of limit linear series have the expected dimension. 
The following definition turns out to determine the behavior of the remaining
components.

\begin{defn}\label{def:weakly-genl} In Situation \ref{sit:binary}, we
say that $X$ is \textbf{weakly Brill-Noether general} for a given $r>0$
if for all $(d_1,d_2)$ with
$$0 \leq d_i \leq g-1, \quad \text{ for } i=1,2$$
the space of $\fg^r_{(d_1,d_2)}$s on $X$ has the expected dimension
$\rho$.
\end{defn}

The reason for our terminology is that the range of multidegrees we consider 
is smaller than would be required in order to obtain a proof of the 
Brill-Noether theorem considering only naive linear series, as is done
by Caporaso in \cite{ca3}. However, our main result on binary curves is
that using our theory of limit linear series, the above multidegrees
are sufficient.

\begin{cor}\label{cor:binary} In Situation \ref{sit:binary}, any component
of a limit linear series space on $X$ for which the general member has
nonconstant multivanishing sequences has the expected dimension $\rho$.

Moreover, if for a given $r$ we have that $X$ and all its partial 
normalizations are weakly Brill-Noether general, then for the same $r$, all 
spaces of limit linear series on $X$ have the expected dimension $\rho$.
\end{cor}

\begin{proof} According to Theorem \ref{thm:exp-dim}, it is enough to
fix all discrete invariants, including choices of multivanishing sequences,
and to show that the gluing conditions impose the maximal codimension. We
simplify notation as follows: for $i=1,2$ let $d_i$ denote the $d_v$ of the
proof of Theorem \ref{thm:exp-dim} applied to the components $Z_i$, and
similarly let $D^i_{\bullet}$ be the sequences of effective divisors on 
$Z_i$, indexed from $0$ to $b+1$. Set $\sL^i:=\sO_{Z_i}(d_i)$, and let 
$a^i$ be the multivanishing sequence on $Z_i$ for $i=1,2$. Finally, let the 
sequences $r^i$ be obtained from the $a^i$ as in Definition 
\ref{defn:strong-genl}. It is clearly enough to show that the gluing 
condition (II) of Definition \ref{defn:lls} imposes the correct codimension 
for each choice of line bundle gluing, so we fix such a gluing.

Let $G$ be the product over all critical $j$ of
$$G(r^1_j, \sL^1(-D^1_j)/\sL^1(-D^1_{j+1})) \times 
G(r^2_{b-j}, \sL^2(-D^2_{b-j})/\sL^2(-D^2_{b+1-j})),$$
so that our gluing conditions can be expressed as (the preimage of)
a subvariety of $G$.
Now, because the $Z_i$ are rational, for any $j$
with $\deg D^1_{j+1} \leq d_1+1$ we have 
$$H^0(Z_1,\sL^1(-D^1_j))/H^0(Z_1,\sL^1(-D^1_{j+1})) \risom 
\sL^1(-D^1_j)/\sL^1(-D^1_{j+1}).$$
Let $j_1$ be maximal with $\deg D^1_{j} \leq d_1+1$. Then if $j=j_1$,
we have that the space
$H^0(Z_1,\sL^1(-D^1_j))/H^0(Z_1,\sL^1(-D^1_{j+1}))$ imbeds into
$\sL^1(-D^1_j)/\sL^1(-D^1_{j+1})$ as a proper linear subspace. Of course, 
if $j>j_1$, then $H^0(Z_1,\sL^1(-D^1_j))=0$, so we see that in fact
the space of choices of $V^1$ can be expressed as the product over
critical $j$ of $G(r^1_j, W^1_j)$, where 
$W^1_j \subseteq \sL^1(-D^1_j)/\sL^1(-D^1_{j+1})$ is a subspace which
is the entire space for $j<j_1$, and zero for $j>j_1$. Letting $j_2$
be maximal such that $\deg D^2_{j} \leq d_2+1$, we have an analogous
statement for $Z_2$, with subspaces 
$W^2_j \subseteq \sL^2(-D^2_j)/\sL^2(-D^2_{j+1})$. Thus, if $P$ is the
space of choices of $(V^1,V^2)$ with the required multivanishing sequences,
we see that $P$ is in fact a subvariety of $G$ which is itself a product. 
It is thus enough to work one index at a time, verifying that the image of
$P$ in the relevant index meets the gluing condition subvariety with the
correct codimension. 

Now, if we have 
$j$ with $j<j_1$ and $b-j<j_2$,
the image of $P$ is all of
$G(r^1_j, \sL^1(-D^1_j)/\sL^1(-D^1_{j+1})) \times 
G(r^2_{b-j}, \sL^2(-D^2_{b-j})/\sL^2(-D^2_{b+1-j}))$ (and this is 
independent of any conditions on other values of $j$), so there is 
certainly no problem with gluing codimension for such values of $j$. If 
$j=j_1$ and $b-j<j_2$,
the image is
$G(r^1_j, W^1_j) \times G(r^2_{b-j}, \sL^2(-D^2_{b-j})/\sL^2(-D^2_{b+1-j}))$,
and it is straightforward to check that the gluing codimension is still
as desired. 
The same still holds if $j<j_1$ and $b-j=j_2$. Values of $j$ with $j>j_1$
or $b-j>j_2$ cannot occur in the multivanishing sequence, so the only case
left to consider is that $j=j_1$ and $b-j=j_2$. The expected codimension of
gluing in this case is only positive if $\deg D^1_j$ occurs in $a^1$ (in
which case it is necessarily $a^1_r$), and if $\deg D^2_{b-j}$ occurs in
$a^2$ (in which case it is $a^2_r$). 

However, we see that this can only happen under very restrictive 
circumstances: indeed, if $b-j_1=j_2$, then for $j>j_1$ we cannot have
$\deg D^1_{j}$ occuring in $a^1$, while for $j<j_1$ we have $b-j>j_2$,
so $\deg D^2_{b-j}$ cannot occur in $a^2$. On the other hand, 
\eqref{eq:eh-genl} implies that if $a^1_r=\deg D^1_{j_1}$, then
$a^2_0 \geq \deg D^2_{b-j_1}$, so we find that $a^2_{\ell}=\deg D^2_{b-j_1}$
for all $\ell$, and similarly $a^1_{\ell}=\deg D^1_{j_1}$ for all $\ell$.
But in this case, our space of limit linear series (if we now let the gluing
of the line bundle vary again) is almost the same as the space of naive
linear series on $X'$, of multidegree 
$(d_1-\deg D^1_{j_1},d_2-\deg D^2_{b-j_1})$, where $X'$ is the partial
normalization of $X$ at the nodes not occurring in the support of
$\deg D^1_{j_1+1}- \deg D^1_{j_1}$. Indeed, the only difference is that
the latter space does not remember the choice of line bundle gluing at
the normalized nodes, so our limit linear series space is smooth over
the latter space, and the expected codimension of gluing is the same in
both.
Now, because $\deg D^1_{j_1+1} > d_1+1$,
and $\deg D^1_{j_1+1}-\deg D^1_{j_1} \leq g+1$, we have that
$d_1 - \deg D^1_{j_1} \leq g-1$, and similarly 
$d_2 - \deg D^2_{b-j_1} \leq g-1$, so we conclude the desired reduction
statement.
\end{proof}

Putting together Corollary \ref{cor:binary} with a theorem of Caporaso,
we find that we have good behavior for $r=1,2$, as follows.

\begin{cor}\label{cor:binary-2} In the situation of Corollary
\ref{cor:binary}, suppose that $X$ is a general 
binary curve, and $r \leq 2$. Then all moduli spaces of limit linear 
series on $X$ have the expected dimension $\rho$.
\end{cor}

\begin{proof} 
According to Theorem 27 of \cite{ca3}, when $X$ is general, usual linear 
series of multidegrees $(d_1,d_2)$ on $X$ have the expected dimension 
whenever $r\leq 2$ and
$$\frac{d-g-1}{2} \leq d_i \leq \frac{d+g+1}{2}$$
for $i=1,2$. One checks easily that these inequalities are weaker than
$0 \leq d_i \leq g-1$ for $i=1,2$.
The asserted result thus follows from Corollary \ref{cor:binary}.
\end{proof}

\begin{rem}\label{rem:enriched-irrel} The complete independence of 
Corollary \ref{cor:binary-2} from both enriched structures and chain
structures may be a bit surprising, insofar as it differs both from the
behavior of the families considered in \S\S \ref{sec:few-nodes-1}
and \ref{sec:few-nodes-2}
and from the intuition that even a curve in
the closure of a Brill-Noether divisor could behave as a general curve
if approached from a suitably general direction. 

However, we observe that this behavior is in fact somewhat predictable, 
insofar as our theory is a direct generalization of the Eisenbud-Harris 
theory for curves of compact type, and even of usual linear series for 
smooth curves. Indeed, in both cases neither enriched structures nor chain 
structures are relevant, and yet non-general curves need not be
Brill-Noether general.
\end{rem}

\begin{ex}\label{ex:binary-hyperelliptic} 
While we hope that the condition $r \leq 2$ in Corollary \ref{cor:binary-2}
may be removed via further analysis, the generality condition on $X$ is
certainly necessary. Indeed, if $X$ is obtained by gluing together two
copies of the same marked rational curve, we see that $X$ always has a
$\fg^1_2$ of multidegree $(1,1)$, which may be considered as a limit linear
series with $b=0$. As long as $g>2$, we have $\rho<0$, so this violates
Corollary \ref{cor:binary-2}. Note that this happens regardless of 
enriched structure or chain structure, so again underlines the peculiar 
irrelevance of enriched structures in the case of binary curves, discussed 
in Remark \ref{rem:enriched-irrel}.
\end{ex}

\begin{rem}\label{rem:binary-explicit} 
While the way we have presented Corollary \ref{cor:binary-2} requires
an imprecise generality hypothesis, in fact one can analyze the situation
quite explicitly, since the weak Brill-Noether generality hypothesis
amounts to studying whether we can find $(r+1)$-dimensional spaces of
polynomials of degree $d_1$ and $d_2$ which can be made to glue to one
another at the chosen points, under suitable choices of gluings. This
can be set up as an explicit linear algebra problem, leading to the 
potential for more precise criteria for Brill-Noether generality,
especially in fixed $r$. For instance, one can easily check via this 
approach 
that the existence of a $g^1_{(1,1)}$ (the only 
possibility for a $g^1_2$) is equivalent to the marked points on $Z_1$ and 
$Z_2$ being the same, up to automorphism.
\end{rem}

\begin{rem}\label{rem:high-d} Note that despite the inequalities in the
definition of being weakly Brill-Noether general, it does not follow from
Corollary \ref{cor:binary} that all spaces of limit linear series on 
binary curves have dimension $\rho$ when $d>2g-2$. The reason is that
the reduction process in the proof of Corollary \ref{cor:binary} may
decrease the degree, or put differently, one can imbed spaces of 
lower-degree linear series into spaces of higher-degree limit linear series
by imposing extra multivanishing conditions. This is not new to our 
notion of limit linear series -- in fact, exactly the same phenomenon 
occurs for Eisenbud-Harris limit linear series.
\end{rem}

\bibliographystyle{amsalpha}
\bibliography{gen}

\end{document}